\documentclass[10pt,twoside]{article}
\usepackage{graphicx}
\usepackage{amsmath}
\usepackage{Latex-document}
\usepackage{latexsym}
\usepackage{amssymb}

 \def\NN{{\mathbb N}}

 \def\ZZ{{\mathbb Z}}

\def\La{\Lambda}

\def\Om{\Omega}

  \def\cH{{\cal H}}  
\def\cC{{\cal C}}   \def\cO{{\cal O}} \def\cU{{\cal U}}
   \def\cP{{\cal P}} \def\cV{{\cal V}}
\def\cE{{\cal E}}    \def\cW{{\cal W}}
\def\cF{{\cal F}}   \def\cR{{\cal R}}

\newtheorem{theo}{Theorem}
\newtheorem{theor}{Theorem}

\newtheorem{corol}{Corollary}
\newtheorem{defi}{Definition}

\newtheorem{conj}{Conjecture}

\newtheorem{prob}{Problem}

\newtheorem{rema}{Remark}

\title{\bf {\boldmath $C^1$}-Generic Dynamics:\vskip -2mm
 Tame and Wild Behaviour\vskip 6mm}

\author{C. Bonatti\vspace*{-0.5cm}\thanks{
Laboratoire de Topologie  UMR 5584 du CNRS, Universit\'e de
Bourgogne, B.P. 47 870, 21078 Dijon Cedex, France. E-mail:
bonatti@u-bourgogne.fr}}
\date{\vspace{-8mm}}

\begin{document}
\maketitle

\thispagestyle{first} \setcounter{page}{265}

\begin{abstract}\vskip 3mm
This paper gives a survey of recent results on the maximal transitive sets of $C^1$-generic diffeomorphisms.

\vskip 4.5mm


\end{abstract}

\vskip 12mm

\section{Introduction} \label{section 1}\setzero
\vskip-5mm \hspace{5mm}

In order to give a global description of a dynamical system (diffeomorphism or flow) on a compact manifold $M$, the first step consists in characterizing the parts of $M$ which are, in some sense, indecomposable for the dynamics. This kind of description will be much more satisfactory if these indecomposable sets are finitely many, disjoint, isolated, and not fragile (that is, persistent in some sense under perturbations of the dynamics).

For non-chaotic dynamics, this role can be played by  periodic orbits, or by minimal sets. However many (chaotic) dynamical systems have infinitely many periodic orbits and a uncountable number of minimal sets. In order to structure the global dynamics using a smaller number of (larger) sets, we need to relax the notion of indecomposability. A weaker natural notion of topological/dynamical indecomposability is the notion of transitivity.

\subsection{Maximal and saturated transitive sets~}
\vskip-5mm \hspace{5mm}

An invariant compact set $K$ of a diffeomorphism $f$ is {\em transitive} if there is a point in $K$ whose positive orbit is dense in $K$. An equivalent definition is the following: for any open subsets $U,V$ of $K$ there is $n>0$ such that $f^n(U)\cap V\neq \emptyset$.

One easily verifies that the closure of the union of an increasing family of transitive sets is a transitive set. Then Zorn's Lemma implies that any transitive set is contained in a {\em maximal transitive set} (i.e. maximal for the inclusion).

However, maximal transitive sets are not necessarily disjoint. For this reason, we also consider the stronger notion of {\em saturated transitive sets}: a transitive set $K$ is saturated if any transitive set intersecting $K$ is contained in $K$. So two saturated transitive sets are always equal or disjoint.

\vskip 2mm
These notions are motivated by  Smale's approach for hyperbolic dynamics, and more specifically for his spectral decomposition theorem (see \cite{Sm}):

\subsection{Smale's spectral decomposition theorem}
\vskip-5mm \hspace{5mm}

For an Axiom A diffeomorphism $f$ on a compact manifold $M$, the set $\Om(f)$ of the non-wandering points is the union of finitely many compact disjoint (maximal and saturated) transitive sets $\La_i$, called {\em the basic pieces}, which are uniformly hyperbolic.

If furthermore $f$ as no cycles, there is a {\em filtration $\emptyset=M_{k+1}\subset M_k\subset\cdots\subset M_1=M$ adapted to $f$}, that is: the $M_i$ are submanifolds with boundary, having the same dimension as $M$, and are strictly  $f$-invariant: $f(M_i)$ is contained in the interior $\stackrel{o}M_i$ of  $M_i$.
Moreover  $\La_i$ is the {\em maximal invariant set} of $f$ in  $M_i\setminus \stackrel{o}M_{i+1}$, that is  $\La_i=\bigcap_{n\in \ZZ} f^n(M_i\setminus \stackrel{o}M_{i+1})$.

Finally this presentation is {\em robust}: the same filtration remains adapted to  any diffeomorphism $g$ in a $C^1$-neighborhood of $f$, and the maximal invariant sets $\La_i(g) =\bigcap_{n\in \ZZ} g^n(M_i\setminus \stackrel{o}M_{i+1})$ are the basic pieces of $g$.

\subsection{A global picture of {\boldmath $C^1$}-generic diffeomorphisms}
\vskip-5mm \hspace{5mm}

It is known from the sixties (see \cite{AbSm,Si}) that Axiom A diffeomorphisms are not $C^1$-dense in $Diff^1(M)$ if $dim(M)>2$ and, of course,  general dynamical systems do not admit such a nice global presentation of their dynamics. One would like to give an analogous description for  as large as possible a class of diffeomorphisms.

In this paper I will present a collection of works, trying to give a coherent global picture of the dynamics of $C^1$-generic diffeomorphisms. There are two types of generic behaviours:
\begin{description}
\item{-} either the manifold contains infinitely many regions having independent dynamical behaviours (we will speak of a {\em wild diffeomorphism}, and in Section~\ref{s.wild} will give  examples of such
behaviours),
\item{-} or one has a description of the dynamics identical to those given by the spectral decomposition theorem: we speak of a {\em tame diffeomorphism}. In this case the role of basic pieces is played by the {\em homoclinic classes} (see the definition in Section~\ref{ss.homoclinic}).
\end{description}

An important class of examples of tame dynamics are the robustly transitive dynamics and in Section~\ref{ss.example} we summarize the known examples of robustly transitive dynamics. The basic pieces of a tame diffeomorphism present a weak form of hyperbolicity called {\em dominated splitting} and {\em volume partial hyperbolicity} (see Section~\ref{ss.hyperbolic} and \ref{ss.bdp}). In Section~\ref{ss.description} we  try to summarize the dynamical consequences of the dominated splittings.

\section{Tame and wild dynamics} \label{section 2}
\setzero \vskip -5mm \hspace{5mm}

\subsection{{\boldmath $C^1$}-generic diffeomorphisms}

\vskip -5mm \hspace{5mm}

In this paper,  we will consider the set of diffeomorphisms $Diff^1(M)$ endowed with the $C^1$-topology. The choice of the topology comes from the fact that most of the perturbating results (Pugh's closing Lemma \cite{Pu}, Hayashi's connecting Lemma and its generalizations \cite{Ha,Ha2,Ar,WX}) are only known in this topology.

Recall that a property $\cP$ is {\em generic} if it is verified on a residual subset $\cR$ of $Diff^1(M)$ (i.e.  $\cR$ contains the intersection of a countable family of dense open subsets).  In this work  we will often use a practical abuse of language; we say:

\centerline{ ``{\em Any $C^1$-generic diffeomorphism verifies
$\cP$}" }

\noindent instead of:

\centerline{``{\em There is a residual subset $\cR$ of $Diff^1(M)$
such that any $f\in\cR$ verifies $\cP$.}" } \vskip 2mm

Let me first recall a famous and classical example, relating Pugh's closing lemma to generic dynamics:
\begin{theor}{\rm \cite{Pu}}\hskip 4mm Let  $f$ be a diffeomorphism on a  compact manifold and $x\in \Om(f)$ be a non-wandering point.  There is $g$ arbitrarily $C^1$-close  to $f$ such that $x$ is periodic for $g$.
\end{theor}

Using a Kupka-Smale argument (genericity of hyperbolicity of the periodic points and the transversality of invariant manifolds) one get:

\begin{corol}
The non-wandering set  $\Om(f)$ of a generic diffeomorphism $f$ is the closure of the set of periodic points of $f$, which are all hyperbolic.
\end{corol}

\subsection{Homoclinic classes}\label{ss.homoclinic}
\vskip -5mm \hspace{5mm}

Let $f$ be a diffeomorphism on a compact manifold and $p$ be a hyperbolic periodic point of $f$ of saddle type. Let $W^s(p)$ and $W^u(p)$ denote the invariant manifold of the orbit of $p$. The {\em homoclinic class} $H(p,f)$ of $p$ is by definition the closure of the transverse intersection points of its invariant manifold:
$$
H(p,f)=\overline{W^s(p,f)\pitchfork W^u(p,f)}.
$$

The homoclinic class $H(p,f)$ is a transitive set canonically associated to the orbit of the periodic point $p$.

There is an other way to see the homoclinic class of $p$: we tell that a periodic point $q$ of saddle type and of same Morse index as $p$ is {\em homoclinically related to $p$} if $W^u(q)$ cuts transversally $W^s(p)$ in at least one point and reciprocally $W^s(q)$ cuts transversally $W^u(p)$ in at least one point. The $\lambda$-lemma (see \cite{Pa}) implies that this relation is an equivalence relation and $H(p,f)$ is the closure of the set of periodic points homoclinically related to $p$.

For Axiom A diffeomorphisms, the homoclinic classes are precisely the basic pieces of Smale's spectral decomposition theorem. However, one easily build examples of diffeomorphisms whose homoclinic classes are not maximal transitive sets. Moreover, B. Santoro \cite{Sa} recently build examples of diffeomorphisms on a $3$-manifold having periodic points whose homoclinic classes are neither disjoint nor equal.

\subsection{Homoclinic classes of generic dynamics}
\vskip-5mm \hspace{5mm}

Conjectured during a long time, Hayashi's connecting lemma allowed the control the perturbations of the invariant manifolds of the periodic points, opening the door for the understanding of generic dynamics.

\begin{theo} {\rm \cite{Ha}}\hskip 4mm Let $p$ and $q$ be two hyperbolic  periodic points of some diffeomorphism $f$. Assume that there is a sequence $x_n$ of points  converging  to a point $x\in W^u_{loc}(p)$  and positive iterates $y_n=f^{m(n)}(x_n)$, $m(n)\geq 0$, converging to a  point $y\in W^s_{loc}(q)$.

Then there is $g$, arbitrarily $C^1$-close to $f$, such that  $x$ and $y$ belong to a same heteroclinic orbit of $p$ and $q$; in other words:

 \centerline{$x\in W^u_{loc}(p,g)$,  $y\in W^s_{loc}(q,g)$ and there is  $n>0$ such that $g^n(x)=y$.}
\label{t.Hayashi}
\end{theo}

If the periodic points $p$ and $q$ in Theorem~\ref{t.Hayashi} belong to a same transitive set, then the sequences $x_n$ and $y_n$ are given by a dense orbit.
In \cite{BoDi2}, using in an essential way Hayashi connecting lemma, we proved:
\begin{theo} For any $C^1$-generic diffeomorphism, two periodic orbits belong to the same transitive set if and only if their homoclinic classes coincide.
\label{t.cmp}
\end{theo}

Motivated by this result  we conjectured:

\vskip 2mm
\noindent
 {\em The homoclinic classes of a generic diffeomorphism coincide with its maximal transitive sets.}

\vskip 2mm
We know now that this conjecture, as stated above, is wrong: in \cite{BoDi3} (see Section~\ref{s.wild}) we show that any manifold $M$ with dimension $>2$ admits a non-empty $C^1$-open subset $U\subset Diff^1(M)$ on which generic diffeomorphisms have an uncoutable family of maximal (an saturated) transitive sets without periodic points.

However, one part of the conjecture is now proved.
Generalizations of Hayashi Connecting Lemma (see \cite{Ha2,Ar,WX})  recently allowed to show:
\begin{theo}  For any $C^1$-generic diffeomorphism, the homoclinic class of any periodic point is a maximal (see \cite{Ar}) and saturated (see \cite{CMP}) transitive set.
\end{theo}

The proof of this theorem is decomposed in two main steps: first,
\cite{Ar} shows that for a generic diffeomorphism $f$ the
homoclinic class of a point $p$ coincides with the intersection
of the closure of its invariant manifolds:
$$
H(p,f) \stackrel{def}{=}\overline{W^s(p,f)\pitchfork W^u(p,f)} =\overline{W^s(p,f)}\cap\overline{W^u(p,f)}.
$$
Then  \cite{CMP} shows that for a generic diffeomorphism $f$ the closure $\overline{W^u(p,f)}$ is Lyapunov stable (and so admits a base of invariant neighborhoods) and $\overline{W^s(p,f)}$ is Lyapunov stable for $f^{-1}$. As a consequence a dense orbit of  a transitive set $T$ intersecting $H(p,f)$ is capted in arbitrarilly small neighborhoods of $\overline{W^u(p,f)}$ and of $\overline{W^u(p,f)}$, proving that $T$ is contained in  $\overline{W^s(p,f)}\cap\overline{W^u(p,f)}$, finishing the proof of the theorem.
\vskip 2mm

\subsection{Tame and wild diffeomorphisms}
\vskip-5mm \hspace{5mm}

Using Theorem~\ref{t.cmp} and the fact that the homoclinic class $H(p,f)$ of a periodic point varies lower semi-continuously with $f$, \cite{Ab} shows the existence of a $C^1$-residual subset $\cR$ of diffeomorphisms (or flows),  such that the cardinality of the set of homoclinic classes is locally constant on $\cR$: for any Kupka-Smale diffeomorphism $f$ let $n(f)\in \NN\cup\{\infty\}$ denote the cardinal of the set of different homoclinic classes $H(p,f)$ where $p$ is an hyperbolic periodic point of $f$; then  any $f\in\cR$ has a $C^1$-neighborhood $U_f$ such that  any $g\in\cR\cap U_f$ verifies $n(g)=n(f)$.

This result induces a natural dichotomy the residual set $\cR$:
\begin{description}
\item{-} a diffeomorphism $f\in \cR$ is {\em tame} if it has finitely many homoclinic classes.
\item{-} a diffeomorphism $f\in \cR$ is {\em wild} if it has infinitely many homoclinic classes.
\end{description}

\section{Tame dynamics}
\vskip-5mm \hspace{5mm}

\subsection{Filtrations, robust transitivity and generic transitivity }
\vskip-5mm \hspace{5mm}

\cite {Ab} shows that the global dynamics of tame diffeomorphisms admit   a good reduction to the dynamics of the transitive pieces (up to reduce the residual set $\cR$). Let $f\in\cR$ be a tame diffeomorphism, then :
\begin{enumerate}
\item
as in the Axiom A case, the non-wandering set is the union of
finitely many disjoint homoclinic classes $H(p_i,f)$.
\item there is a filtration $\emptyset=M_{k+1}\subset M_k\subset\cdots\subset M_1=M$ adapted to $f$ such that $H(p_i,f)$ is the maximal invariant set in $M_i\setminus \stackrel{o}M_{i+1}$.
\item\label{i.generic} moreover (up to reduce the open neighborhood $U_f$ defined above) this filtration holds for any $g\in U_f$, and for $g\in \cR\cap U_f$ the maximal invariant set of $g$ in $M_i\setminus \stackrel{o}M_{i+1}$ is the homoclinic class $H(p_{i,g},g)$.
\item there is a good notion of attractors: either a homoclinic class is a topological attractor (that is, its local basin contains a neighborhood of it) or its stable manifold has empty interior. Then the union of the basin of the attractors of $f$ is a dense open set of $M$ (see \cite{CaMo}).
\end{enumerate}

The item~\ref{i.generic} above shows that the transive sets $H(p_i,f)$ are not fragil. In a previous work, \cite{DPU} introduced the following notion:
\begin{defi}\label{d.robust} Let $f$ be a diffeomorphism of some compact manifold $M$. Assume that there is some open set $U\subset M $ and a $C^1$-neighborhood $\cV$ of $f$ such that, for any $g\in \cV$, the maximal invariant set $\La_g=\bigcap_{n\in\ZZ} g^n(\bar U)$ is a compact transitive set contained in $U$.

Then $\La_f$ is called a {\em robustly transitive set} of $f$.
\end{defi}

In the definition above, if one has $U=M$ (and so $\La_g=M$ for
any $g\in\cV$), then $f$ is called a {\em robustly transitive
diffeomorphism}.

This notion is slightly stronger that the property given by the item~\ref{i.generic} above; so we have to relax  Definition~\ref{d.robust}: we say that $\La_f$ is {\em generically transitive} if, in the notations of Definition~\ref{d.robust}, the maximal invariant set  $\La_g$ is transitif for $g$ in a residual subset of $\cV$.

At this moment, there are no known examples of generic transitive sets which are not robustly transitive. So it is natural to ask if this two notions are equivalent:

\centerline{\em Generic transitivity
$\stackrel{?}{\Longleftrightarrow}$ robust transitivity ?}

\subsection{Examples of robust transitivity}\label{ss.example}
\vskip-5mm \hspace{5mm}

The Axiom A dynamics are obvious examples of tame dynamics. On compact surfaces, tame diffeomorphisms are, in fact, Axiom A diffeomorphisms, but there are many non-hyperbolic examples in higher dimensions.

Even if this talk is mostly devoted to diffeomorphisms, let us observe that the most famous robustly transitive non-hyperbolic attractor is the Lorenz attractor (geometric model, see \cite{GuWi,ABS}) for flows on $3$-manifolds. There are many generalizations of this attractor, called singular attractors, for flow on $3$-manifolds, see for instance \cite{MPP}. See also \cite{BPV} for robust singular attractors in dimension greater ou equal than $4$, having a singular point with Morse index (dimension of unstable manifold) greater than $2$.

The first example of non-Anosov robustly transitive diffeomorphism is due to Shub \cite{Sh}: it is a diffeomorphisms on the torus $T^4$ which is a skew product over an Ansov map on the torus $T^2$, such that the dynamics on the fibers is dominated by the dynamics on the basis.

Then Ma\~n\'e \cite{Man2} built an example  of robustly transitive non-hyperbolic diffeomorphism on the torus $T^3$ by considering a bifurcation of an Anosov map  $A$ having $3$ real positive different eigenvalues $\lambda_1<1<\lambda_2<\lambda_3$ : he performs a saddle node bifuraction creating a (new) hyperbolic saddle of index $1$ (breaking the hyperbolicity) in the weak unstable manifold of a fixed point (of index $2$) of the Anosov map $A$.

Then \cite{BoDi} shows that a diffeomorphism $f$ admits $C^1$ perturbations which are robustly transitive, if $f$ is:
\begin{enumerate}
\item  the time one diffeomorphism of  any transitive Anosov flow.
\item  the product $(A,id)$ where $A$ is some Anosov map and $id$ is the identity map of any compact manifold.
\end{enumerate}

The second case can be easily generalized to any skew product of
an Anosov map by rotations of the circle $S^1$. The same
technique also allows \cite{BoDi} to build example of robustly
transitive attractors, by perturbating product maps of any
hyperbolic attractors by the identity map of some compact
manifold.

Each of these previous example was partially hyperbolic (see the definitions in Section~\ref{ss.hyperbolic}): they admits a splitting $TM=E^s\oplus E^c\oplus E^u$ where $E^s$ is uniformly contracting and $E^u$ is uniformly expanding, and it was  conjectured that partial hyperbolicity was a necessary condition for robust transitivity. Then \cite{BoVi} generalizes Ma\~n\'e example above and exhibits robustly transitive difeomorphisms  on $T^3$ having a uniformly contracting $1$-dimensional bundle, but no expanding bundle (there is a splitting $TM=E^s\oplus E^{cu}$), and robustly transitive difeomorphisms on $T^4$ having no hyperbolic subbundles (neither expanding nor contracting): there just admits  an invariant dominated splitting $TM=E^{cs}\oplus E^{cu}$).

We do not known what are the manifolds admitting robustly transitive diffeomorphism. For instance:
\begin{conj}There is no robuslty transitive diffeomorphism on the sphere
$S^3$.
\end{conj}

This conjecture has been proved in \cite{DPU} assuming the
existence of a codimension $1$ (center stable or center unstable)
foliation, using Novikov Theorem. Notice that all the known
examples of robustly transitive diffeomorphisms on $3$-manifolds
admits an invariant codimension $1$ foliation. However this
conjecture remains still open.

\subsection{Dominated splitting and partial hyperbolicity: definitions}\label{ss.hyperbolic}
\vskip-5mm \hspace{5mm}

Let $f$ be a $C^1$-diffeomorphism of a compact manifold and let
$\cE$ be an $f$-invariant compact subset of  $M$. Let  $TM_x=
E_1(x)\oplus\cdots\oplus E_k(x)$, $x\in \cE$, be a splitting of
the tangent space at any point of   $\cE$. This splitting is a
{\em dominated splitting} if it verifies the following properties:
\begin{enumerate}
\item For any $i\in\{1,\dots,k\}$, the dimension $E_i(x)$ is independent of
$x\in\cE$.
\item The splitting is $f_*$-invariant (where $f_*$ denots the differential of $f$):  $E_i(f(x))= f_*(E_i(x))$.
\item There is $\ell\in \NN$ such that , for any $x\in \cE$, for any  $1\leq i<j\leq k$ and any  $u\in E_i(x)\setminus\{0\}, v\in E_j(x)\setminus\{0\}$  one has:
$$
\frac{\|f_*^\ell(u)\|}{\|u\|}\leq\frac{\|f_*^\ell(v)\|}{2\|v\|}.
$$
\end{enumerate}
\begin{rema}
\begin{description}{\rm
\item{-} (Continuity) Any dominated splitting on a set $\cE$ is continuous and extend in a unique way to the closure $\bar\cE$.
\item{-} (Extension to a neighborhood) There is a neighborhood $U$ of $\bar\cE$ on which the maximal invariant set  $\Lambda(\bar U,f)$ has a dominated splitting extending those on $\cE$.
\item{-} (Robust) There is a $C^1$-neighborhood  $\cU_f$ of $f$ such that, for any $g\in \cU_f$, the maximal invariant set   $\La(\bar U,g)$ has a dominated splitting varying continuously with $g$.
\item{-} (Unicity) If $\cE$ has a dominated splitting, then there is a (unique) dominated spliting $TM|_{\cE}= E_1\oplus\cdots\oplus E_k$, called {\em the finest dominated splitting}, such that any other dominated splitting  $F_1\oplus\cdots \oplus F_l$ over  $\cE$ is obtained by grouping the $E_i$ in
packages.}
\end{description}
\end{rema}

One of the  $E_i$ is {\em uniformly contracting} if  (up to increase  $\ell$ in the definition above)
$\frac{\|f_*^\ell(u)\|}{\|u\|}\leq \frac 12$ for all $x\in \cE$ and all $u\in E_i(x)\setminus \{0\}$.
In the same way $E_i$ is {\em uniformly expanding} if $\frac{\|f_*^\ell(u)\|}{\|u\|}\geq \frac 12$ for all $x\in \cE$ and all $u\in E_i(x)\setminus \{0\}$.

An  $f$-invariant compact set $K$ is {\em hyperbolic} if it has a dominated splitting $TM|_K= E^s\oplus E^u$ where $E^s$ is uniformly contracting and $E^u$ is uniformly expanding.  The compact $f$-invariant set  $K$ is {\em partially hyperbolic} if it has a dominated splitting and if at least one of the bundles $E_i$ of its finest dominated splitting  is uniformly contracting or expanding. Let $E^s$ and $E^u$ be the sum of the uniformly contracting and expanding subbundles, respectively, and let $E^c$  be the sum of the other subbundles. One get a new  dominated splitting  $E^s\oplus E^c$, $E^c\oplus E^u$ or $E^s\oplus E^c\oplus E^u$, and these bundles are called the stable, central et unstable bundles, respectively.

An $f$-invariant compact set $K$ is called {\em volume hyperbolic} if there is a dominated splitting whose extremal bundles $E_1$ and $E_k$ contracts and expands uniformly the volume, respectively. Notice if one of these bundle has dimension $1$, it is uniformly contracting or expanding. In particular, a volume hyperbolic set in dimension $2$ is a uniformly hyperbolic set, and in dimension $3$ it is partially hyperbolic (having at least one uniformly hyperbolic bundle).

\subsection{Volume hyperbolicity for the robust transitivity} \label{ss.bdp}
\vskip-5mm \hspace{5mm}

Generalizing previous results by Ma\~n\'e \cite{Man} (in dimension $2$) and by \cite{DPU} in dimension $3$, \cite{BDP} (for robustly transitive set) and \cite{Ab} for generically transitve sets show:
\begin{theo} Any robustly (or generically) transitive set is volume hyperbolic.
\label{t.volume}
\end{theo}

Then any robustly transitive set in dimension $2$ is a hyperbolic
basic set (result of Ma\~n\'e) and in dimension $3$ is partially
hyperbolic (\cite{DPU}). In higher dimension,  the dominated
splitting may have all the subbundles of dimension greater than
$2$, so the expansion or contraction of the volume does no more
imply the hyperbolicity of the bundle, see the example in
\cite{BoVi}.

The proof of Theorem~\ref{t.volume} has  two very different steps
(as in \cite{Man}). The first one consists in showing that the
lake of dominated splitting allows to ``mixe" the eigenvalues of
the periodic orbits, creating an homothecy; a periodic orbit
whose  differential at the period is an homothecy is (up to a
small perturbation) a sink or a source, breaking the
transitivity. For that we just perturb the linear cocycle defined
by the differential of $f$, and then we use a Lemma of Franks
(\cite{Fr}) for realizing the linear perturbation as a dynamical
perturbation. Let state precisely this result:

\begin{theo}{\rm \cite{BDP}}\hskip 4mm Let  $f$ be a diffeomorphism of a compact manifold $M$, and let $p$ be a hyperbolic periodic saddle. Assume that the homoclinic class $H(p,f)$ do not have any dominated splitting . Then, given any $\varepsilon>0$, there is a  periodic point $x$  homoclinically related to $p$, with the following property:

Given any neighborhood $U$ of the orbit of $x$, there is a diffeomorphism $g$,  $\varepsilon$-$C^1$-close to $f$, coinciding with $f$ out of $U$ along the orbit of $x$, such that the  differential $g^n_*(x)$ is a  homothecy, where  $n$ is the periode of $x$.
\label{t.BDP}
\end{theo}

The second step consists in proving the uniform contraction and expansion of the volume in the extremal bundles. As in \cite{Man}, one uses  Ma\~n\'e's Ergodic Closing Lemma  to realize  a lake of uniform expansion (or uniform expansion) of the volume in the extremal bundle by a periodic orbit $z$ of a $C^1$-perturbation of $f$: if furthermore, the differential of this point restricted to the corresponding extremal bundle is an homothecy (as in Theorem~\ref{t.BDP}) one get a sink or a source, breaking the transitivity.

\vskip 5mm
For flows, the existence of singular point lies to additional difficulties. In dimension $3$, \cite{MPP2} show that a robustly transitive set $K$ of a flow on a compact $3$-manifold is a uniformly hyperbolic set if it does not contain any singular point. If $K$ contains a singular point then all the singular points in $K$ have the same Morse index and $K$ is a singular  attractor if this index is $1$ and a singular repellor if this index is $2$ (see also \cite{CMP2}).

\subsection{Topological description of the dynamics with dominated splittings}\label{ss.description}
\vskip-5mm \hspace{5mm}

The dynamics of diffeomorphisms admitting dominated splitting is
already very far to be understood.

In dimension $2$, Pujals and Sambarino (see \cite{PuSa,PuSa2}) give a very precise description of $C^2$-diffeomorphism whose non-wandering set admits a dominated splitting.
\begin{description}
\item{-} the periods of the non-hyperbolic periodic points is upper bounded.
\item{-} $\Om(f)$ is the union of finitely many normally hyperbolic circles on which a power of $f$ is a rotation, (maybe infinitely many) periodic points contained in a finite family of periodic normally hyperbolic segments and finitely many pairwise disjoint homoclinic classes, each of them containing at most finitely many non-hyperbolic periodic orbits.
\end{description}

This result is close to Ma\~n\'e 's result, in dimension $1$, for
$C^2$-maps far from critical points (see \cite{Ma}). We hope that
this result can be generalized in any dimension, for dynamics
having a codimension $1$ strong stable bundle:

\begin{conj} Let $f$ be a $C^2$-diffeomorphism and $K$ be a compact locally maximal invariant set of $f$ admitting a dominated splitting $TM|_K=E^s\oplus F$ where $F$ has dimension $1$ and $E^s$ is uniformly contracting.

Then $K$ is the union of finitely many normally hyperbolic circles on which a power of $f$ is a rotation, of periodic points contained in a union of finitely many  normally hyperbolic periodic intervals and finitely many pairwize disjoint homoclinic classes each of them containing at most finitely many non-hyperbolic periodic points.
\end{conj}

In this direction S. Crovisier \cite{Cr} obtained some progress
in the case where there is a unique non-hyperbolic periodic point.

General dominated splitting cannot avoid wild dynamics: multiplying any diffeomorphism by a uniform contraction and a uniform expansion, we get a normally hyperbolic and partially hyperbolic set. However a dominated splitting give some information of the possible bifurcations and on the index of the periodic point: see \cite{BDPR} which investigate in this direction. In particular a diffeomorphism cannot present any homoclinic tangency if it admits a dominated splitting whose non-hyperbolic bundles are all of dimension $1$. We hope that this kind of dominated splitting avoid wild behaviours, but this is unknown, even in dimension $3$:

\begin{conj} Let $M$ be a compact $3$-manifold and denote by $\cP\cH(M)$ the $C^1$-open set of partially hyperbolic diffeomorphisms of $M$ admitting a dominated splitting $E^s\oplus E^c\oplus E^u$ where all the bundles have dimension $1$.

The open set $\cP\cH(M)$ does not contain any wild diffeomorphism: in other word any generic diffeomorphism in $\cP\cH(M)$ is tame.
\end{conj}

For partially hyperbolic diffeomorphism (having a splittin
$E^s\oplus E^c\oplus E^u$),  Brin and Pesin (\cite{BrPe}) show
the existence of unique foliations $\cF^s$ and $\cF^{u}$,
$f$-invariant and tangent to  $E^s$ and $E^u$ respectively. The
dynamics of the  strong stable and  the strong unstable
foliations play an important role for the understanding of the
topological and ergodical properties of a partially hyperbolic
diffeomorphisms. Let mention two results on these foliations:
\cite{DoWi} shows that a dense open subset of partially
hyperbolic diffeomorphisms (having strong stable and strong
unstable foliaitons) verify the ``accessibility property", that
is, any two points can be joined by a concatenation of  pathes
tangent successively to the strong  stable or the strong unstable
foliations. When the center direction has dimension $1$,
\cite{BDU} shows the minimality of at least one of the strong
stable or strong unstable foliations for a dense open subset of
the robustly transitive systems in $\cP\cH(M)$, where $M$ is a
compact $3$-manifold.

However there is no general result on the existence of invariant foliations tangent to the central bundle even if it has dimension $1$.
When a partially hyperbolic diffeomorphism presents an invariant foliation  $\cF^c$ tangent to the center bundle $E^c$ and which is {\em plaque expansive}, \cite{HPS} shows that this foliations is structurally stable: any $g$ close to $f$ admits a foliation $\cF^c_g$ topologically conjugated to $f$ and such that (up to this conjugacy of foliation) $g$ is isotopic to $f$ along the center-leaves. This gives a very strong rigidity of the dynamics. This deep result was a key step for the construction of the examples of robustly transitive examples in \cite{Sh,Man2,BoDi} (there is now new proofs which do not use the stabilitity of the center foliation (see\cite{Bo})). So an important problem is:

\begin{prob}{\rm  (1)
Does it exist robustly transitive partially hyperbolic
diffeomorphisms having an invariant center foliation which is not
plaque expansive?

(2) If a transitive partially hyperbolic diffeomorphism admits an
invariant center foliation, is it {\em dynamically coherent}?
that is, does it admit  invariant center-stable and center
unstable foliations which intersect along the center foliation?

(3) If the center bundle  is $1$-dimensional, is there an
invariant center foliation?}

\end{prob}

\section{Wild dynamics}\label{s.wild}\setzero

\vskip-5mm \hspace{5mm}

Very little is known on wild diffeomorphisms: for surfaces, it is not known whether $C^1$-wild diffeomorphisms
exist (recall that the Newhouse phenomenon is a $C^2$-generic phenomenon, see \cite{N}).

In dimension $\geq 3$, the known examples are all of them due to the existence of homoclinic classes which do not admit, in a persistent way, any {\em dominated splitting} (see the first examples in \cite{BoDi2}).  Then following the same ideas, \cite{BoDi3} present wild diffeomorphisms exhibiting, in a locally generic way, infinitely many hyperbolic and non-hyperbolic non-periodic attractors . The same example will present maximal transitive sets without any periodic orbits. The rest of this section is devoted to a short presentation of these examples:

Consider an open subset $\cU\subset Diff^1(M)$ such that for any $f$ in $\cU$ there is a periodic point $p_f$ depending continuously on $f$ and verifying:

\begin{description}
\item{-} For all $f\in\cU$ the homoclinic class $H(p_f,f)$ contains two periodic points of different Morse indices, and having each of them a complex (non-real) eigenvalue (this eigenvalue is contracting for one  point and expanding for the other).
\item{-} For all $f\in \cU$ there are two periodic points having the same Morse index as  $p_f$ and homoclinically related to $p_f$ such that the jacobian of the derivative of $f$ at the period is strictly greater than one for one off this point and stricly less than one for the other point.
\end{description}

First item means that the homoclinic class $H(p_f,f)$ do not have any dominated splitting, and that this property is robust. So Theorem~\ref{t.BDP} shows that  $H(p_f,f)$ admits periodic  points whose derivative can be perturbated  in order to get an homothecy.  Second item above allows to choose this point having a jacobian (at the period) arbitrarily close to $1$. Then a new pertubation allows to get a periodic point whose derivative at the period is the identity. Considering then perturbations of the identity map, we get:

\begin{theo}{\rm \cite{BoDi3}}\hskip 4mm There is a residual part  $\cR$ of the open set $\cU$ defined above,  such that any $f\in \cR$ admits an infinite family of periodic disks  $D_n$ (let  $t_n$ denote the period), whose orbits are pairwize disjoint, and verifying the universal following property:

Given any $C^1$-open set $\cO$ of diffeomorphisms from the disk $D^3$ to its interior $\stackrel{o}{D^3}$, there is  $n$ such that the restriction of  $f^{t_n}$ to the disk $D_n$ is smoothly conjugated to an element of $\cO$.
\end{theo}

Notice that the set  of diffeomorphisms  $g\colon D^3 \to\stackrel{o}{D^3}$ contains an open subset $\cU_0$ verifying the property of  $\cU$ described above, one get some kind of renormalisation process: there is a residual part of $\cU$ containing  infinitely many periodic disks $D_n$ containing each of them   infinitely many periodic subdisks   themself containing   infinitely many  periodic subdisks and so on...
In that way one build a tree such that each   branch is a sequence  (decreasing for the inclusion), of strictly periodic orbits of disks whose periods go to infinity, and whose radius go to zero. The intersection of this sequence is a  Lyapunov stable (and so saturated) transitif  compact set, conjugated to an {\em adding machine} (see for instance \cite{BS} for this notion)  and so without periodic orbits. The set of the infinite branches of this  tree is uncountable, given the following result :
\begin{theo} {\rm \cite{BoDi3}} Given any compact manifold  $M$ of dimension $\geq 3$, there is an open subset $\cV$ of $Diff^1(M)$ and a residual part $\cW$ of $\cV$, such that any $f\in\cW$ admits an uncountable family of  saturated transitif sets without periodic orbits.\label{t.sauvage}
\end{theo}

\label{lastpage}


\begin{thebibliography}{aa}
\bibitem{Ab} F. Abdenur,  Generic robustness of spectral decompositions, {\em preprint IMPA}, (2001)
\bibitem{AbSm} R. Abraham and S. Smale, Non-genericity of $\omega$-stability, {\em Global Analysis, vol XIV of Proc. Symp. Pure Path. (Berkeley 1968) Amer. Math Soc.,} (1970).
\bibitem{ABS} V.S. Afraimovitch,  V.V. Bykov and L.P. Shil'nikov  On the appearance and structure of the Lorenz attractor {\em  Dokl. Acad. Sci. USSR}, 234,  (1977), 336--339.

\bibitem{Ar}M.-C. Arnaud,  Cr{\'e}ation de connexions en topologie ${C}^1$, {\em Ergod. Th. {$\&$} Dynam. Systems}, {21}, {(2001)},
{339--381}.
\bibitem{Bo} Ch. Bonatti,  Dynamique g\'en\'eriques: hyperbolicit\'e et transitivit\'e, {\em S\'eminaire Bourbaki} n904, Juin 2002.
\bibitem{BoDi} Ch. Bonatti and L.J. D\'\i az, Persistent nonhyperbolic transitive diffeomorphisms. {\em Ann. of Math.}, (2) 143 (1996), no. 2, 357--396.
\bibitem{BoDi2} Ch. Bonatti and L.J. D\'\i az,  Connexions h\'et\'eroclines et g\'en\'ericit\'e d'une infinit\'e de puits ou de sources, {\em Ann. Scient. \'Ec. Norm. Sup.}, $4^e$ s\'erie, t32, (1999), 135--150.
\bibitem{BoDi3} Ch. Bonatti and  L.J. D\'\i az,  On maximal transitive sets of generic diffeomorphisms, {\em preprint} (2001).
\bibitem{BDP} C. Bonatti, L. J. D\'\i az and  E. Pujals, A $C^1-$generic dichotomy for diffeomorphisms: weak forms of hyperbolicity or infinitely many sinks or sources,
to appear at {\em Annals of Math.}
\bibitem{BDPR} C. Bonatti, L. J. D\'\i az, E. Pujals and J. Rocha,  Robust transitivity and heterodimensional cycles, To appear in {\em Asterisque}.
\bibitem{BDU} C. Bonatti, L. J. D\'\i az and  R. Ures,  Minimality of the strong stable and strong unstable foliations for partially hyperbolic diffeomorphisms, to appear in the {\em Publ. Math. de l'inst. Jussieu}.
\bibitem{BPV} Ch. Bonatti, A. Pumari\~no and M.Viana,  Lorenz attractors with arbitrary expanding dimension, {\em C.R. Acad Sci Paris}, 1, 325, Serie I,  (1997), 863--888.
\bibitem{BoVi} Ch. Bonatti et M.Viana,   SRB measures for partially hyperbolic attractors: the
contracting case, {\em Israel Journal of Math.},{\bf 115},  (2000), 157--193.
\bibitem{BS}
J. Buescu and I. Stewart,
 Liapunov stability and adding machines,
{\em Ergodic Th. $\&$ Dyn. Syst.}, {\bf 15}(2), (1995), 271--290.

\bibitem{BrPe} M. Brin et Ya. Pesin , Partially hyperbolic dynamical systems, {\em Izv. Acad. Nauk. SSSR}, 1, (1974), 177--212.

\bibitem{CaMo} C. Carballo, C. Morales,  Homoclinic classes and finitude of attractors for vector
fields on $n$-manifolds, {\em Preprint} (2001).

\bibitem{CMP} C. Carballo, C. Morales, and M.J. Pac\'\i fico,
 Homoclinic classe for $\cC^1$-generic vector fields,  to appear at {\em Erg. The. and Dyn Sys.}
\bibitem{CMP2} C. Carballo, C. Morales and M.J. Pac\'\i fico,  Maximal transitive sets with singularities for generic $C^1$-vector fields, {\em Boll. Soc. Bras. Mat.} 31, n3,  (2000), 287--303.

\bibitem{Cr} S. Crovisier,  Saddle-node bifurcations for hyperbolic sets, to appear at {\em Erg. Theor and Dyn Sys.}

\bibitem{DPU} L. J. D\'\i az, E. Pujals and R. Ures,  Partial hyperbolicity and robust transitivity {\em Acta Mathematica}  vol. 183, (1999), 1--43.
\bibitem{DoWi} D. Dolgopyat and A. Wilkinson ,Stable accessibility is $C^1$-dense, to appear in {\em Ast\'erisque}.
\bibitem{Fr} J. Franks, Necessary conditions for stability of diffeomorphisms, {\em Trans. A.M.S.}, {\bf 158},  (1971), 301-308.
\bibitem{GuWi} J. Guckenheimer and R.F. Williams, Structural stability of Lorenz attractors, {\em Publ. Math. IHES} 50,  (1979), 59--72.
\bibitem{Ha} S. Hayashi, Connecting invariant manifolds and the  solution of the $C^1$-stability and $\Om$-stability conjectures for flows, {\em Ann. of Math.}, {\bf 145},  (1997), 81--137.
\bibitem{Ha2} S. Hayashi, A $C^1$ make or break lemma,{\em Bol. Soc. Bras. Mat.} 31,(2000), 337--350.
\bibitem{HPS}M. Hirsch, C. Pugh, et M. Shub, {\em Invariant manifolds},
Lecture Notes in Math., 583, Springer Verlag, 1977.

\bibitem{Ma} R.  Ma\~n\'e,  Hyperbolicity, sinks and measure in one-dimensional dynamics, {\em Comm. Math. Phys.} 100,  (1985), 495--524.
\bibitem{Man} R.  Ma\~n\'e,  An ergodic closing lemma, {\em Annals of Math.}  vol. 116, (1982), 503--540
\bibitem{Man2} R.  Ma\~n\'e, Contributions to the stability conjecture, {\em  Topology},  17,  (1978), 386--396.
\bibitem{MPP} C. Morales, M.J. Pacifico and E. Pujals, Singular hyperbolic systems, {\em Proc. Amer. Math. Soc.} {127},(1999), {3393--3401}.
\bibitem{MPP2} C. Morales, M.J. Pacifico and E. Pujals, Robust transitive singular sets for $3$-flows are partially hyperbolic attractors and repellers, {\em preprint IMPA} (1999).

\bibitem{N}
S.Newhouse,
 Diffeomorphisms with infinitely many sinks,
{\em Topology}, {\bf 13}, (1974), 9--18.
\bibitem{Pa} J. Palis, On Morse Smale dynamical systems, {\em Topology}, 8, (1969), 385--405.
\bibitem{Pu} C. Pugh, The closing lemma, {\em Amer. Jour. of Math.}, {\bf 89}, (1967), 956--1009.

\bibitem{PuSa} E. Pujals and M. Sambarino, Homoclinic tangencies and hyperbolicity for surface diffeomorphisms. {\em Ann. of Math.}  151 , no. 3, (2000), 961--1023.
\bibitem{PuSa2} E. Pujals and M. Sambarino, The dynamics of dominated splitting, {\em Preprint IMPA} (2001).
\bibitem{Sa} B. Santoro , Colis\~ao, colapso e explos\~ao  de classes holocl\'\i nicas, {\em Thesis PUC Rio de Janeiro} (2001).
\bibitem{Sh}M. Shub, Topological transitive diffeomorphism on $T^4$, {\em Lect.  Notes in Math.},  {\bf 206},   39  (1971).
\bibitem{Si} R. Simon , A $3$-dimensional Abraham-Samle example, {\em Proc. Amer. Math Soc.}, 34, (1972), 629--630.
\bibitem{Sm} S. Smale,  Differentiable dynamical systems, {\em Bull. Am. Math. Soc.}, 73, (1967), 747--817.
\bibitem{WX}{L. Wen and Z. Xia},  $C^1$ connecting lemmas, {\em Trans. Amer. Math. Soc.}, {352}, {(2000)}, {5213--5230}.
\end{thebibliography}
\end{document}